\documentclass[12pt]{amsart}
\usepackage{amssymb,latexsym,eufrak,amsmath,amscd, graphicx}

\usepackage{amsfonts}
\usepackage{amscd}

\usepackage[matrix,arrow]{xy}

\author{Florin Ambro} 
\address{RIMS, Kyoto University\\
Kyoto 606-8502, JAPAN.}
\email{ambro@kurims.kyoto-u.ac.jp}

\newcommand{\Q}{{\mathbb Q}}
\newcommand{\Z}{{\mathbb Z}}

\newcommand{\R}{{\mathbb R}}

\newcommand{\cC}{{\mathcal C}}

\newcommand{\cL}{{\mathcal L}}
\newcommand{\cM}{{\mathcal M}}

\newcommand{\cO}{{\mathcal O}}

\newcommand{\cR}{{\mathcal R}}
\newcommand{\cS}{{\mathcal S}}

\newcommand{\emb}{\operatorname{emb}}
\newcommand{\Ker}{\operatorname{Ker}}

\newcommand{\relint}{\operatorname{relint}}
\newcommand{\Proj}{\operatorname{Proj}}

\theoremstyle{plain}
\newtheorem{thm}{Theorem}[section]
\newtheorem{theom}{Theorem}

\newtheorem{lem}[thm]{Lemma}

\theoremstyle{definition}
\newtheorem{defn}[thm]{Definition}
\newtheorem{exmp}[thm]{Example}

\newtheorem{ack}{Acknowledgments}   

\begin{document}

\bibliographystyle{amsalpha+}
\title[Toric FGA]
{Asymptotically saturated toric algebras}

\begin{abstract}
We show the finite generation of certain invariant 
graded algebras defined on toric weak log Fano fibrations. 
These are the toric version of FGA algebras, 
recently introduced by Shokurov in connections to the 
existence of flips.
\end{abstract}

\maketitle



\section*{Introduction}


\footnotetext[1]{This work was supported by a 
Twenty-First Century COE Kyoto Mathematics Fellowship,
and by the Japan Society for Promotion of Sciences.
}
\footnotetext[2]{1991 Mathematics Subject Classification. 
Primary: 14E30. Secondary: 14M25, 11J25.}

Shokurov's new approach for establishing
the existence of flips is to reduce this problem, by 
induction on dimension, to the finite generation of 
graded algebras which are {\em asymptotically saturated} 
with respect to weak log Fano fibrations (\cite{PLflips}). 
Shokurov showed that these algebras are finitely generated 
in dimension one and two, and conjectured this to be 
true in any dimension (\cite{PLflips}). Our main result 
is the positive answer to the toric case of this conjecture. 

\begin{theom}\label{main}
Let $\pi\colon X\to S$ be a proper surjective toric morphism 
with connected fibers, and let $B$ be an invariant $\Q$-divisor 
on $X$ such that $(X,B)$ has Kawamata log terminal  
singularities and $-(K+B)$ is $\pi$-nef.
\begin{itemize}
\item[(1)] Let 
$
\cL\subseteq \bigoplus_{i=0}^\infty\pi_*\cO_X(iD)
$ 
be an invariant graded $\cO_S$-subalgebra which is
asymptotically saturated with respect to $(X/S,B)$,
where $D$ is an invariant divisor on $X$.
Then $\cL$ is finitely generated.
\item[(2)] The number of rational maps 
$X\dashrightarrow \Proj(\bar{\cL})$, where $\bar{\cL}$ 
is the normalization of an $\cO_S$-algebra as in (1), 
is finite up to isomorphism.
\end{itemize}
\end{theom}

The toric case of asymptotic saturation, the key property 
ensuring finite generation in (1), can be explicitely 
written down as a Diophantine system (see the proof of 
Theorem~\ref{maini}). To see this in a special case,
let $S$ be a point and let $X=T_N\emb(\Delta)$ be 
a torus embedding for some lattice $N$. Let $M$ be 
the lattice dual to $N$, and consider a compact convex set 
of maximal dimension $\square$ in $M_\R$. This defines 
a toric graded algebra
$$
\cR(\square)=\bigoplus_{i=0}^\infty
(\bigoplus_{m\in M\cap i\square}{\mathbb C}\cdot \chi^m),
$$
which is finitely generated if and only if $\square$ is a
rational polytope.
On the other hand, the log discrepancies of the log pair 
$(X,B)$ with respect to toric valuations can be encoded in a 
positive function $\psi\colon N_\R\to \R$, and $\psi$ 
determines a rule to enlarge any convex set in $M_\R$ to
an open convex neighborhood. 
The asymptotic saturation of $\cR(\square)$ with respect 
to $(X,B)$ means that the lattice points of the 
neighborhood of $j\square$ are already contained in 
$j\square$, for every sufficiently divisible positive 
integer $j$. This Diophantine property restricts the way 
that $\square$ can be approximated with rational points 
from the outside.

The key technical tool behind Theorem~\ref{main} 
is a known result in Geometry of Numbers, an effective 
bound on the width of a convex set in terms of the number of 
lattice points it contains (\cite{KL}).

The outline of this paper is as follows. In Section 4 
we explicitely describe toric asymptotic saturation and reduce 
Theorem~\ref{main} to its special case when the algebra 
is normal and associated to a convex set, the equivalent of 
$\cR(\square)$ above. The rest of the paper is devoted to 
this special case. In Section 1 we collect some elementary 
results on convex sets and their support functions, and on 
Diophantine approximation. In Section 2 we characterize 
asymptotic saturation in geometric terms (Theorem~\ref{chr}) 
and obtain a boundedness result (Theorem~\ref{klbound}). 
These are used in Section 3 to prove 
Theorem~\ref{main}, by induction on dimension.

\begin{ack} I would like to thank Nobuaki Sugimine
for useful discussions. 
\end{ack}


\section{Preliminary}


We collect in this section elementary results on convex sets,
toric geometry and Directed Diophantine Approximation, which 
we will need later. We refer the reader to Oda~\cite{Oda88} 
for basic notions and terminology on toric varieties and 
convex sets.

Throughout this section, $N$ is a lattice, with
dual lattice $M$. We have a duality pairing
$\langle \cdot,\cdot \rangle\colon M_\R\times N_\R\to \R$,
defined over $\Z$. 


\subsection{Convex sets and support functions}


Fix a {\em convex rational polyhedral cone} $\sigma\subseteq 
N_\R$, that is $\sigma$ is spanned by finitely many elements of 
$N$. We denote by $\cS(\sigma)$ the set of all functions 
$h\colon \sigma\to \R$ satisfying the following properties:
\begin{itemize}
\item[1)] {\em positively homogeneous}: $h(t e)=t\cdot h(e)$
for $t\ge 0, e\in \sigma$.
\item[2)] {\em upper convex}: $h(e_1+e_2)\ge h(e_1)+
h(e_2)$ for $e_1,e_2\in \sigma$.
\end{itemize}

\begin{thm}\label{Rocka}
\begin{itemize}
\item[(i)] Every function $h\in \cS(\sigma)$ 
is continuous.
\item[(ii)] Let $(h_i)_{i\ge 1}$ be a sequence of 
functions in $\cS(\sigma)$ which converges pointwise, 
and set
$
h(e)=\lim_{i\to \infty} h_i(e)
$
for 
$e\in \sigma$.
Then $h\in \cS(\sigma)$, and the sequence $(h_i)_i$ 
converges uniformly to $h$ on compact subsets of 
$\sigma$.
\end{itemize}
\end{thm}

\begin{proof} This is a special case of~\cite{Rock70},
Theorems 10.1 and 10.8. 
\end{proof}
For a function $h\colon \sigma\to \R$, define
$$
\square_h=\{m\in M_\R; \langle m,e\rangle\ge h(e),
\forall e\in \sigma\}.
$$
A {\em convex polytope} $K\subset M_\R$ is the convex hull 
of a finite set in $M_\R$. A {\em rational convex polytope} 
is the convex hull of a finite set in $M_\Q$.
A {\em rational convex polyhedral set} is the intersection
of finitely many rational affine half spaces in $M_\R$.
We denote by $\cC(\sigma^\vee)$ the set of all nonempty 
closed convex sets $\square\subseteq M_\R$
satisfying the following two properties:
\begin{itemize}
\item[1)] $\square+\sigma^\vee=\square$;
\item[2)] $\square\subseteq K+\sigma^\vee$, for some convex 
polytope $K\subset M_\R$.
\end{itemize}
The {\em support function} $h_\square\colon \sigma\to \R$
of $\square\in \cC(\sigma^\vee)$ is defined by
$$
h_\square(e)=\inf_{m\in \square}\langle m,e\rangle.
$$
 
\begin{thm}\label{dictionary} 
The maps $\square\mapsto h_\square$ and $h\mapsto \square_h$ 
are inverse to each other, inducing a bijection
$
\cC(\sigma^\vee)\simeq \cS(\sigma).
$
Under this correspondence, the Minkowski sum $\square+\square'$ 
and a nonnegative scalar multiple $t\square$ correspond to 
$h_\square+h_{\square'}$ and $t h_\square$, respectively.
\end{thm}

 We omit the proof of this theorem, being similar to that of 
~\cite{Oda88}, Theorem A.18. When $\sigma=N_\R$, this is the usual 
correspondence between compact convex sets and support functions.
Note that $K+\sigma^\vee\in \cC(\sigma^\vee)$, for every 
compact convex set $K\subset M_\R$, but not all elements of 
$\cC(\sigma^\vee)$ are of this form. Such an example is
$\sigma^\vee=\{(x,y)\in \R^2; x,y\ge 0\}$ and 
$\square=\{(x,y)\in \sigma^\vee; xy\ge 1\}$. 
Nevertheless, we have

\begin{lem}\label{rp} The following properties are equivalent for 
$\square\in \cC(\sigma^\vee)$:
\begin{itemize}
\item[(i)] $\square$ is a rational convex polyhedral set.
\item[(ii)] $\square=K+\sigma^\vee$ for some rational convex 
polytope $K$ such that no vertex of $K$ belongs to the 
Minkowski sum of $\sigma^\vee$ and the convex hull of the 
other vertices of $K$.
\end{itemize}
Furthermore, assume $\dim(\sigma)=\dim(N)$. Then $K$ is 
uniquely determined by $\square$.
\end{lem}

\begin{exmp} If $\sigma=N_\R$, then $K=\square$. If
$\sigma$ is the positive cone in $\R^d$ and $\square$ 
is a Newton polytope, then $K$ is the convex 
hull of the compact faces of $\square$.
\end{exmp}


\subsection{Proper toric morphisms with affine base}

Toric morphisms with affine base which are
proper, surjective and with connected fibers, are in one to
one correspondence with fans having convex support.

Indeed, let $\Delta$ be a fan in a lattice $N$ such that its
support $\vert \Delta\vert=\bigcup_{\tau\in \Delta}
\tau$ is a convex rational polyhedral cone. Let 
$\bar{N}=N/(N\cap \vert \Delta\vert\cap(-\vert\Delta\vert))$
and let $\bar{\sigma}\subset \bar{N}_\R$ be the image of 
$\vert \Delta\vert$ under the natural projection. Then 
$
T_N\emb(\Delta)\to T_{\bar{N}}(\bar{\sigma})
$
is a toric morphism with affine base, which is proper, 
surjective, with connected fibers. 

Conversely, let $\pi\colon X\to S$ be a proper toric 
morphism of toric varieties, with $S$ affine. Thus 
$X=T_N\emb(\Delta)$, $S=T_{N'}(\sigma')$ and $\pi$ 
corresponds to a lattice homomorphism 
$\varphi_\Z\colon N\to N'$ such that $\Delta$ is a 
finite fan in $N$, $\sigma'$ is a strongly convex 
rational polyhedral cone in $N'$ and 
$\vert\Delta\vert=\varphi^{-1}(\sigma')$. 
In particular, $\vert \Delta\vert$ is a convex rational 
polyhedral cone.
Then $\varphi_\Z$ factors through $\bar{N}=
N/(N\cap \vert \Delta\vert\cap(-\vert\Delta\vert))$ and we 
have a commutative diagram
\[ \xymatrix{
X \ar[dr]_\pi \ar[rr] & & T_{\bar{N}}(\bar{\sigma})
\ar[dl]^j \\
          & S   &
} \] 
where $j$ is finite on its image.

Let now $D=\sum_{e\in \Delta(1)} d_e V(e)$ be an invariant 
$\Q$-divisor on $X$ which is $\Q$-Cartier. This means that 
there exists a function $h\colon \vert\Delta\vert\to \R$ 
such that $h$ is $\Delta$-linear and $h(e)=-d_e$ for every 
$e\in \Delta(1)$. In particular, $h$ is positively homogeneous.
The $\Q$-divisor $D$ is $\pi$-nef if 
$h$ is upper convex; it is $\pi$-ample if  
for every maximal cone $\sigma\in \Delta$ there 
exists $m_\sigma\in \square_h$ such that
$\sigma=\{e\in \vert\Delta\vert; h(e)=\langle 
m_\sigma,e\rangle\}$.


\subsection{The ample fan of a convex rational 
polyhedral set}


To each convex rational polyhedral set $\square\subseteq M_\R$
we associate a fan $\Delta_\square$ in a quotient lattice of $N$,
as follows. Assume first that $\dim(\square)=\dim(M)$. Let $K$ 
be a rational polytope associated to $\square$ by Lemma~\ref{rp}, 
with vertices $v_1,\ldots,v_l$. The support function of $\square$ 
is 
$h(e)=\min_{j=1}^l \langle v_j,e\rangle$, and the cones
$$
\{e\in \vert\Delta\vert; \langle v_j,e\rangle=h(e)\}
\ (1\le j\le l)
$$
are the maximal dimensional cones of a fan $\Delta_\square$ 
in $N$. The support $\vert\Delta_\square\vert$ is the unique 
convex cone $\sigma\subseteq N_\R$ such that 
$\square\in \cC(\sigma^\vee)$.

If $\dim(\square)<\dim(M)$, choose a point
$m_0\in M_\Q\cap \square$ and denote $\square'=\square-m_0$.
Define
$$
N'=N/(N\cap {\square'}^\perp).
$$
If $M'$ is the dual lattice of $N'$, then $M'_\R$ can be
identified with the smallest vector subspace of $M_\R$
which contains $\square'$. We have $\dim(\square')=\dim(M')$
and let $\Delta_{\square'}$ be the fan in $N'$ defined as 
above. This fan is independent of the choice of 
$m_0$, and we denote it again by $\Delta_\square$. Its support
is a convex set.

\begin{defn} $\Delta_\square$ is called the {\em ample fan} 
of the rational convex polyhedral set $\square\subset M_\R$.
\end{defn}

Assume now that $\pi\colon X\to S$ is a proper toric morphism
with affine base $S$. We may write $X=T_N\emb(\Delta)$ and 
$S=T_{\bar{N}}(\bar{\sigma})$, and $\pi$ corresponds to a 
lattice homomorphism $\varphi_\Z\colon N\to \bar{N}$ such that
$\vert\Delta\vert=\varphi^{-1}(\bar{\sigma})$.

For a rational polyhedral convex set 
$\square\in \cC(\vert\Delta\vert^\vee)$, the $T_N$-invariant 
$\cO_S$-algebra
$$
\cR(\square)=\bigoplus_{i=0}^\infty
(\bigoplus_{m\in M\cap i\square}{\mathbb C}\cdot \chi^m).
$$
is normal and finitely generated. The induced toric rational
map
\[ \xymatrix{
X \ar[dr] \ar@{.>}[rr]^\Phi & & \Proj(\cR(\square))\ar[dl] \\
          & S   &
} \] 
is defined over $S$, and $\Proj(\cR(\square))$ is the torus 
embedding of the ample fan $\Delta_\square$.
If $\dim(\square)=\dim(M)$, then $\Delta_\square$ is a fan in $N$
with $\vert \Delta_\square\vert=\vert\Delta\vert$, hence $\Phi$ is 
birational in this case. The invariant
$\Q$-divisor $\sum_{e\in \Delta_\square(1)}-h(e)V(e)$ is ample 
relative to $S$. If $\dim(\square)<\dim(M)$, $\Delta_\square$ is a 
fan in $N'$, whose support is the image of $\vert \Delta\vert$ under
the natural projection.


\subsection{Directed Diophantine Approximation}

Let $m\in M_\R$ and let $e\in N_\R$.
Let $I$ be a positive integer and let $\Vert\cdot \Vert$
be a norm on $M_\R$.

\begin{thm}[cf.~\cite{BSh03}]\label{appusual} 
For every $\epsilon>0$, there exists 
a positive multiple $k$ of $I$ and there exists $\bar{m}\in M$ 
such that $\langle \bar{m}-km,e\rangle\in (-\epsilon,0]$ and 
$\Vert \bar{m}-km\Vert<\epsilon$.
\end{thm}

\begin{proof} We may find a decomposition $M=M'\oplus M''$,
with dual decomposition $N=N'\oplus N''$, such that
$m=m'+m''$, $m'\in M'_\Q, m''\in M''_\R$ and 
$\{e''\in N''; \langle m'',e''\rangle\in \Q\}=\{0\}$.
Let $k_1$ be a positive integer such that $I\vert k_1$ and 
$k_1m'\in M'$. Since 
$$
\{e''\in N''; \langle k_1m'',e''\rangle\in \Q\}=\{0\},
$$ 
we infer by \cite{Cassels57}, Chapter III, Theorem IV, that
the subgroup generated by the class of $k_1m''$ is
dense in the torus $M''_\R/M''$. Equivalently, the
set $\bigcup_{j\ge 1}(M''+jk_1m'')$ is dense in $M''_\R$.
In particular, the following system has a solution for 
some $j\ge 1$:
\[
\left\{ \begin{array}{l}
 m''_j\in M'' \\
 \langle m''_j+j k_1 m'',e''\rangle = \langle m''_j+j k_1 m'',
 e\rangle \in [0,\epsilon)\\
 \Vert m''_j+j k_1 m''\Vert<\epsilon.
\end{array} 
\right. \]
Then $k=j k_1$ and $\bar{m}=km'-m''_j$ satisfy the
desired properties.
\end{proof}

\begin{thm}\label{app1} 
Assume that 
$e\notin \{e'\in N; \langle m,e'\rangle\in \Q\}\otimes_\Z \R$.
Then for every $\epsilon>0$, there exists a 
positive multiple $k$ of $I$ and there exists 
$\bar{m}\in M$ such that $\langle \bar{m}-km,e\rangle
\in (-\epsilon,0)$ and $\Vert \bar{m}-km\Vert<\epsilon$.
\end{thm}

\begin{proof} We may decompose $M=M'\oplus M''$, with dual
decomposition $N=N'\oplus N''$, such that $m=m'+m''$,
$m'\in M'_\Q, m''\in M''_\R$ and $\{e''\in N''; \langle m'',
e''\rangle\in \Q\}=\{0\}$. In particular,
$$
N'_\R=\{e'\in N; \langle m,e'\rangle\in \Q\}\otimes_\Z \R.
$$
Let $e=e'+e''$ be the unique 
decomposition with $e'\in N'_\R$ and $e''\in N''_\R$. 
Our assumption means that $e''\ne 0$. Let $k_1$ be a 
positive integer such that $I\vert k_1$ and $k_1m'\in M'$. 
We have
$$
\{e''\in N''_\Q; \langle k_1m'',e''\rangle \in \Q\}=\{0\}.
$$
By~\cite{Cassels57}, Chapter III, Theorem IV, 
the subgroup generated by the class of $k_1m''$ is
dense in the torus $M''_\R/M''$. Equivalently, the
set $\bigcup_{j\ge 1}(M''+jk_1m'')$ is dense in $M''_\R$.
Since $e''\ne 0$, 
$$
\{m''\in M''_\R; \langle m'',e''\rangle \in (0,\epsilon),
\Vert m''\Vert<\epsilon\}
$$
is a non-empty open subset of $M''_\R$. Therefore the 
following system has a solution for some $j\ge 1$:
\[
\left\{ \begin{array}{l}
 \bar{m}_j''\in M'' \\
 \langle \bar{m}_j''+j k_1 m'',e''\rangle\in (0,\epsilon)\\
 \Vert \bar{m}_j''+ j k_1 m''\Vert<\epsilon.
\end{array} 
\right. \]
The claim holds for $k=j k_1$ and $\bar{m}=km'-\bar{m}_j''$.
\end{proof}

\begin{lem}\label{app2} 
If $\langle m,e\rangle\in \Q$, the following properties are 
equivalent:
\begin{itemize}
\item[(i)]
$m\in \{m'\in M;\langle m',e\rangle \in \Q\}\otimes_\Z \R$.
\item[(ii)]
$e\in \{e'\in N;\langle m,e'\rangle \in \Q\}\otimes_\Z \R$.
\end{itemize}
\end{lem}

\begin{proof} Assume that (i) holds.
We may find a decomposition $N=N'\oplus N''$, with dual 
decomposition $M=M'\oplus M''$, such that $e=e'+e''$,
$e'\in N'_\Q,e''\in N''_\R$ and 
$$
\{m''\in M'';\langle m'',e''\rangle \in \Q\}=\{0\}.
$$
The assumption means that $m\in M'_\R$. We may find a 
decomposition $M'=M'_1\oplus M'_2$, with dual decomposition
$N'=N'_1\oplus N'_2$, such that $m=m'_1+m'_2$,
$m'_1\in M'_{1,\Q}$, $m'_2\in M'_{2,\R}$ and
$$
\{e'_2\in N'_2;\langle m'_2,e'_2\rangle \in \Q\}=\{0\}.
$$
We have
$$
\langle m'_2,e'\rangle=
\langle m,e'\rangle-\langle m'_1,e'\rangle=
\langle m,e\rangle-\langle m'_1,e'\rangle\in \Q.
$$
Therefore $e'\in N'_{1,\Q}$. Let $e''=\sum_i r_i e''_i$,
where $r_i\in \R$ and $\{e''_i\}_i$ is a basis of $N''$. 
Then $e',e''_i\in N_\Q$, 
$\langle m,e'\rangle \in \Q$ and $\langle m,e''_i\rangle=0$.
Therefore
$$
e=e'+e''\in \{e'\in N;\langle m,e'\rangle \in \Q\}\otimes_\Z \R,
$$
that is (i) holds. The statement is symmetric in $m$ and $e$, 
hence the converse holds as well.
\end{proof}


\section{Asymptotic saturation}


Throughout this section, we fix a lattice $N$ and a 
convex rational polyhedral cone $\sigma \subseteq N_\R$.

\begin{defn} A {\em log discrepancy function} 
is a function 
$
\psi\colon \sigma \to \R
$ 
satisfying the following properties:
\begin{itemize}
\item[(i)] $\psi$ is positively homogeneous.
\item[(ii)] $\psi(e)>0$ for $e\ne 0$.
\item[(iii)] $\psi$ is continuous.
\end{itemize}
\end{defn}

\begin{exmp}\label{gld} Let $\Delta$ be a fan in $N$ with
$\vert \Delta \vert=\sigma$. Let 
$B=\sum_{e\in \Delta(1)} b_e V(e)$ be an invariant  
$\R$-divisor on $X=T_N\emb(\Delta)$ such that $K+B$ 
is $\R$-Cartier and the pair $(X,B)$ has Kawamata 
log terminal singularities. Equivalently, there 
exists a function
$
\psi\colon \sigma \to \R
$
such that $\psi(e)=1-b_e>0$ for every 
$e\in \Delta(1)$, and $\psi$ is $\Delta$-linear.
Then $\psi$ is a {\em log discrepancy function}.

The terminology comes from the following property:
let $e\in N^{prim}\cap \sigma$ be a primitive lattice 
point, corresponding to a toric valuation $v_e$ of $X$. 
Then $\psi(e)$ is the log discrepancy of $(X,B)$ at $v_e$.
\end{exmp}

\begin{lem}\label{4th} Let $\psi\colon \sigma\to \R$ be a 
log discrepancy function. Then the set
$
\{e\in \sigma;\psi(e)\le 1\}
$
is compact.
\end{lem}

\begin{proof} Choose a norm $\Vert \cdot \Vert$ on
$N_\R$.
Since $\psi$ is a log discrepancy function,
the infimum
$
c_0=\inf\{\psi(e); e\in \sigma, \Vert e\Vert=1\}
$
is a well defined positive real number. We have
$
\psi(e)\ge c_0\Vert e\Vert,
$ 
for $e\in \sigma$. In particular, 
$$
\{e\in \sigma; \psi(e)\le 1\}\subseteq
\{e\in \sigma; \Vert e\Vert\le c_0^{-1}\}.
$$
The left hand side is a closed set, since $\psi$ is continuous,
and the right hand side is a bounded set. Therefore the claim
holds.
\end{proof}

\begin{defn} For an arbitrary function $h\colon \sigma\to \R$, 
define
$$
\stackrel{\circ}{\square}_h=
\{m\in M_\R; \langle m,e\rangle> h(e),\forall e\in 
\sigma\setminus 0\}.
$$
\end{defn}

\begin{defn} Let $\square\in \cC(\sigma^\vee)$ and let
$\psi\colon \sigma\to \R$ be a log discrepancy function.
We say that 
\begin{itemize}
\item[-] $\square$ is {\em $\psi$-saturated} if
$
M\cap \stackrel{\circ}{\square}_{h_\square-\psi}\subset 
\square,
$
where $h_\square\in \cS(\sigma)$ is the support function
of $\square$. Note that $\square=\square_{h_\square}$.
\item[-] $\square$ is {\em asymptotically $\psi$-saturated}
if there exists a positive integer $I$ such that $j\square$ 
is $\psi$-saturated, for every $I\vert j$.
\end{itemize}
\end{defn}

Note that saturation (asymptotic saturation) is invariant 
under lattice (rational) translations of the convex set.

\begin{thm}[Characterization of asymptotic saturation]\label{chr} 
Let $\square\in \cC(\sigma^\vee)$ be a rational polyhedral 
set and let $\psi\colon \sigma\to \R$ be a log discrepancy 
function. 

Let $N''=\{e\in N; \square\ni m \mapsto \langle m,e\rangle
\in \R \mbox{ constant}\}$, with dual lattice $M''$, 
and define $\psi''\colon N''_\R\to \R$ by
$$
\psi''(e)=\psi(e).
$$ 
The ample fan $\Delta_\square$ is a fan in $N'=N/N''$ 
with support $\sigma'=\pi(\sigma)$, where 
$\pi\colon N_\R\to N'_\R$ is the natural projection. Define
$\psi'\colon \sigma'\to \R$ by 
$$
\psi'(e')=\inf_{e\in \pi^{-1}(e')}\psi(e).
$$
Then $\square$ is asymptotically $\psi$-saturated if and only 
if the following hold:
\begin{itemize}
\item[(1)] $M''\cap \stackrel{\circ}{\square}_{-\psi''}=\{0\}$.
\item[(2)] $\Delta_\square(1)\subset \{e'\in N'_\R;
\psi'(e')\le 1\}$.
\end{itemize}
\end{thm}

\begin{proof} After a rational translation, we may assume
$0\in \square$. In particular, $N''=N\cap \square^\perp$ 
and $h(e)=h'(\pi(e))$, where $h\in \cS(\sigma)$ and 
$h'\in \cS(\sigma')$ are the support functions of 
$\square\subset M_\R$ and $\square\subset M'_\R$, respectively.

Assume that (1) and (2) hold. Fix a positive integer $j$ 
such that $I\vert j$ and $jh(N)\subseteq \Z$ and assume 
that $m\in M$ satisfies
$$
\langle m,e\rangle>(jh-\psi)(e), \forall e\in \sigma\setminus 0.
$$
Choose a decomposition $M=M'\oplus M''$, and decompose $m=m'+m''$. 
Since $h\vert_{N''}=0$, we obtain $m''\in M''\cap Q''$, hence 
$m''=0$ by (1). In particular, we have
$$
\langle m,e'\rangle>jh'(e')-\psi'(e'), 
\forall e'\in \sigma'\setminus 0.
$$
For every $e'\in \Delta_\square(1)$, we have
$\langle m,e'\rangle\in \Z$, hence (2) gives
$\langle m,e'\rangle\ge jh'(e')$. Since $h'$
is $\Delta_\square$-linear, we obtain
$$
\langle m,e'\rangle\ge jh'(e'), 
\forall e'\in \vert\Delta_\square\vert,
$$
hence $m\in \square_{jh'}$. 
Therefore $jh$ is $\psi$-saturated.

For the converse, assume that $jh$ is $\psi$-saturated
for every $I\vert j$. We first check (1). 
Fix $m''\in M''\cap Q''$. Let $\Vert\cdot \Vert$ be a norm
on $M_\R$ which is compatible with the decomposition
$M=M'\oplus M''$.
Since $\psi$ is continuous, there exists $\epsilon>0$ 
such that
$$
\langle m'',e\rangle+\psi(e)>0, \forall e\in 
S(\sigma), \Vert e'\Vert < \epsilon.
$$
The rational convex polyhedral set $\square$ has the same 
dimension as $M'_\R$. Therefore there exists 
$m'\in M'\cap \relint(k\square)$, for some positive 
integer $k$. We have
$$
\langle m',e'\rangle>k h'(e'), \forall e'\in \sigma'\setminus 0.
$$ 
The continuity of $\psi$ implies that the following number
is well defined:
$$
t=-\inf\{\frac{\langle m'',e\rangle+\psi(e)}
{\langle m',e'\rangle-kh'(e')}; e\in S(\sigma), 
\Vert e'\Vert\ge \epsilon\}.
$$
Let $j$ be a positive multiple of $I$ such that $j>t$. 
The identity
$$
\langle jm'+m'',e\rangle-(jkh-\psi)(e)=
j(\langle m',e'\rangle-kh'(e'))+\langle m'',e\rangle+\psi(e)
$$
implies that
$$
\langle jm'+m'',e\rangle>(jkh-\psi)(e), \forall e\in 
S(\sigma).
$$
Since $jkh$ is $\psi$-saturated, we infer that
$jm'+m''\in \square_{jkh}$. In particular, $jm'+m''\in M'$, 
hence $m''=0$. This proves (1).

For (2), fix $e'\in \Delta_\square(1)$ and assume by
contradiction that $\psi'(e')> 1$. We may find a 
basis $e_1,\ldots,e_d$ of $N$ with $e_1=e'$. Let
$\Vert \cdot \Vert$ be the absolute value norm on
$N_\R$ with respect to this basis and denote
$$
S(\sigma)=\{e\in \sigma; \Vert e\Vert=1\}.
$$
The face 
$
\{m\in \square; \langle m,e_1\rangle=h(e_1)\}
$
of $\square$ is a positive dimensional convex 
polyhedral set, hence there exists a $1$-dimensional 
rational compact convex set $\square_1$ with 
$$
\square_1\subset \relint(\{m\in \square; \langle 
m,e_1\rangle=h(e_1)\}).
$$
It is easy to see that there exists a positive real 
number $t_1$ such that 
$$
M\cap t\square_1\ne \emptyset \mbox{ for }t> t_1.
$$ 
Consider the following set
$$
C=\{e\in S(\sigma); \psi(e)\le \langle e_1^*,e\rangle\}.
$$
Since $\psi$ is continuous, the (possibly empty) set 
$C$ is closed. Furthermore, $e_1\notin C$ and $\square_1$
is included in the relative interior of the face of 
$\square$ corresponding to $e_1$, hence 
$\langle m,e\rangle-h(e)>0$ for $e\in C$ and $m\in \square_1$. 
We infer that the following number is well defined
$$
t_2=\sup_{m\in \square_1,e\in C}
\frac{(e_1^*-\psi)(e)}{\langle m,e\rangle-h(e)}.
$$
Let $j$ be a positive multiple of $I$ such that 
$j>\max(t_1,t_2)$. Since $j>t_1$, there exists $m_j\in M$
such that $m_j+e_1^*\in j\square_1$. We have 
$$
\langle m_j,e\rangle-(j h-\psi)(e)=
j(\langle \frac{m_j+e_1^*}{j},e\rangle-h(e))
-(e_1^*-\psi)(e).
$$
Since $j>t_2$, we obtain $m_j\in \stackrel{\circ}{\square}_{jh-\psi}$.
Since $j\square$ is $\psi$-saturated, we obtain
$m_j\in \square_{jh}$. This is a contradiction,
since
$$
\langle m_j,e_1\rangle=jh(e_1)-1<jh(e_1).
$$
This proves (2).
\end{proof}

\begin{thm}\label{klbound}
Let $\psi\colon N_\R\to \R$ be a log discrepancy 
function such that $-\psi$ is upper convex and 
$M\cap \stackrel{\circ}{\square}_{-\psi}=\{0\}$.
Then there exists $e\in N\setminus 0$ such that 
$$
\psi(e)+\psi(-e)\le C,
$$ 
where $C$ is a positive constant depending only on $\dim(N)$.
\end{thm}

\begin{proof} Let $\square=\square_{-\frac{\psi}{2}}$. Since
$\psi$ is positive, we have $0\in \square\subset 
\stackrel{\circ}{\square}_{-\psi}$.
Then $\square$ is a compact convex set, of dimension 
$\dim(N)=d$, with support function $-\frac{\psi}{2}$, such 
that $M\cap \square=\{0\}$. 
By~\cite{KL}, Theorem 4.1, there exists $e\in N\setminus 0$
such that 
$$
\max_{m\in \square}\langle m,e\rangle-
\min_{m\in \square}\langle m,e\rangle\le
c_0 d^2 \lceil \sqrt[d]{1+\#(M\cap \square)}\rceil,
$$ 
where $c_0$ is a positive universal constant and 
$\#(M\cap \square)$ is the number of lattice points of
$\square$. In our case, this means
$$
\psi(e)+\psi(-e)\le C=2c_0 d^2 \lceil \sqrt[d]{2}\rceil.
$$
\end{proof}

\begin{thm}[Toric Asymptotic CCS]\label{asyccs} 
Let $\psi\colon \sigma \to \R$ be a log discrepancy 
function. We denote by $\cM(\psi)$ the set of
rational polyhedral sets $\square\in \cC(\sigma^\vee)$
such that 
\begin{itemize}
\item[(1)] $\square$ is asymptotically $\psi$-saturated.
\item[(2)] $h_\square-\psi$ is upper convex.
\end{itemize}
Then the set of ample fans $\{\Delta_\square\}_{
\square\in \cM(\psi)}$ is finite.
\end{thm}

\begin{proof} 
Let $\square\in \cM(\psi)$, with support function $h$. 
After a rational translation, we may assume $0\in \square$.
Let $N''=N\cap \square^\perp$ and let $d=\dim(N'')$.
If $d=0$, that is $\dim(\square)=\dim(M)$, the ample 
fan $\Delta_\square$ is a fan in $N$ with 
$\vert\Delta_\square\vert=\sigma$, and by 
Theorem~\ref{chr} we have 
$$
\Delta_\square(1)\subseteq 
N^{prim}\cap \{e\in N_\R; \psi(e)\le 1\}.
$$
The right hand side is a finite set, hence we infer that
the number of fans $\Delta_\square$ is finite.

Assume now $d>0$. We will show that $N''$ belongs to
a finite set of sublattices of $N$. Since $h\vert_{N''}=0$, 
$
-\psi\vert_{N''_\R}=(h-\psi)\vert_{N''_\R}
$
is an upper convex function.
By assumption, 
$M''\cap \stackrel{\circ}{\square}_{-\psi\vert_{N''_\R}}
=\{0\}$.
By Theorem~\ref{klbound}, there exists $e_1\in N''\setminus 0$
such that
$
\psi(e_1)+\psi(-e_1)\le C.
$
We may assume that $e_1$ is a primitive element of $N$.
Consider the lattice $N'=N/(\Z\cdot e_1)$ and let 
$\pi_\Z\colon N\to N'$ be the induced projection map. 
There exists $h'\colon \sigma'\to \R$ such that $h=h'\circ \pi$. 
Define
$$
\psi'(e')=\inf_{e\in \pi^{-1}(e')}\psi(e).
$$
Then $\psi'$ is a log discrepancy function on $N'_\R$
and $\square=\square_{h'}\in \cM(\psi')$, by 
Lemma~\ref{restriction}. We repeat this argument $d$ times, 
until we obtain a basis $e_1,\ldots,e_d$ of $N''$ with the 
following properties:
\begin{itemize}
\item[(i)] $\psi(e_1)+\psi(-e_1)\le C$.
\item[(ii)]
$
\inf(\psi\vert_{e_k+\sum_{i=1}^{k-1}\R e_i})+
\inf(\psi\vert_{-e_k+\sum_{i=1}^{k-1}\R e_i})\le C,
$
for $2\le k\le d$.
\end{itemize}
By Lemma~\ref{4th}, $e_1$ belongs to a finite set.
By Lemmas~\ref{restriction} and~\ref{4th}, 
$e_k$ belongs to a finite set modulo $\sum_{i=1}^{k-1}\Z e_i$,
for every $k$. 
Therefore $N''$, the subspace of $N$ generated by $e_1,\ldots,e_d$,
belongs to a finite set of sublattices of $N$.

For $N''$ as above, let $N'=N/N''$. There exists 
$h'\colon \sigma'\to \R$ such that $h=h'\circ \pi$. Define
the log discrepancy function $\psi'\colon \sigma'\to \R$ by
$$
\psi'(e')=\inf_{e\in \pi^{-1}(e')}\psi(e).
$$
By Theorem~\ref{chr} again, we have
$
\Delta_\square(1)\subseteq {N'}^{prim}\cap\{e'\in N'_\R;
\psi'(e')\le 1\}.		
$
Therefore the ample fans $\Delta_\square$ are finitely many.

Since $d\le \dim(\sigma)$, we conclude that the number of 
ample fans is finite.
\end{proof}

\begin{lem}[Restriction of saturation]\label{restriction} 
Let $\psi\colon \sigma\to \R$ be a log discrepancy function,
let $\square\in \cC(\sigma^\vee)$ and let $\pi_\Z\colon N\to N'$ 
be a quotient lattice. We identity the dual lattice $M'$ with 
$M\cap \Ker(\pi)^{\perp} \subset M$. 
The image $\sigma'=\pi(\sigma)$ is a rational convex polyhedral 
cone in $N'_\R$.

Assume that $m_0\in M\cap \square$. The convex set 
$\square'=(\square-m_0)\cap M'_\R$ belongs to $\cC({\sigma'}^\vee)$ 
and its support function $h'\in \cS(\sigma')$ is computed as
follows
$$
h'(e')=\sup\{h(e)-\langle m_0,e\rangle; 
e\in \sigma\cap\pi^{-1}(e')\}.
$$
Define a positively homogeneous function 
$\psi'\colon \sigma'\to \R$ by
$$
\psi'(e')=h'(e')-\sup\{h(e)-\langle m_0,e\rangle-\psi(e)
;e\in \sigma\cap\pi^{-1}(e')\}.
$$
Then the following properties hold:
\begin{itemize}
\item[(i)] If $\square$ is $\psi$-saturated, 
then $\square'$ is $\psi'$-saturated.
\item[(ii)] For a positive integer $k$, define
$\psi'_k\colon \sigma'\to \R$ by
$$
\psi'_k(e')=kh'(e')-\sup\{kh(e)-\langle m_0,e\rangle-\psi(e)
;e\in \sigma\cap\pi^{-1}(e')\}.
$$
If $\square$ is asymptotically $\psi$-saturated, then 
$\square'$ is asymptotically $\psi'_k$-saturated.
\item[(iii)] If $h-\psi$ is upper convex, then 
$h'-\psi'$ is upper convex and $\psi'$ is a log discrepancy 
function.
\item[(iv)] $\psi'(e')\ge \inf_{e\in \sigma\cap\pi^{-1}(e')}\psi(e)$.
\end{itemize}
\end{lem}

\begin{proof} We may assume $m_0=0$ after a translation of $\square$.

(i) The inclusion 
$
\stackrel{\circ}{\square}_{h'-\psi'}\subseteq 
\stackrel{\circ}{\square}
$
is easy to see. Since $\square$ is $\psi$-saturated, 
$M'\cap \stackrel{\circ}{\square}_{h-\psi}\subset \square$. 
Therefore $M'\cap \stackrel{\circ}{\square}_{h'-\psi'}
\subseteq M'\cap \square=\square'$.

(ii) Note first the following identity:
$$
(\psi'_j-\psi'_k)(e')=(j-k)\sup_{\pi(e)=e'} h(e)+
\sup_{\pi(e)=e'} (kh-\psi)(e)-
\sup_{\pi(e)=e'} (jh-\psi)(e).
$$
Therefore $\psi'_k\le \psi'_j$ for $k\le j$.

By assumption, there exists a positive integer $I$
such that $jh$ is $\psi$-saturated for every $I\vert j$.
Fix $k\ge 1$ and let $j$ be a common multiple of $I$ and
$k$. By (i), $jh'$ is $\psi'_j$-saturated. Since
$\psi'_j\ge \psi'_k$, we infer that $jh'$ is also
$\psi'_k$-saturated.

(iii) The upper convexity of $h'-\psi'$ follows from the
upper convexity of $h-\psi$ and the formula
$$
(h'-\psi')(e')=\sup_{\pi(e)=e'}(h-\psi)(e).
$$
In particular, $\psi'$ is a continuous function, being
the diference of the continuous functions
$h'$ and $h'-\psi'$. Furthermore, $\psi'$ is positively 
homogenous by its definition.
Let $0\ne e'\in \sigma'$. The 
restriction $\psi\vert_{\pi^{-1}(e')}$ is strictly 
positive, continuous and at least $1$ outside some 
bounded subset, by Lemma~\ref{4th}. Therefore
$$
\inf_{\pi(e)=e'}\psi(e)>0.
$$
We conclude from (iv) that $\psi'(e')>0$.

(iv) This is a direct consequence of the definitions
of $h'$ and $\psi'$.
\end{proof}


\section{Rational polyhedral criterion}


\begin{thm}\label{1rational} 
Let $\sigma\subseteq N_\R$ be a rational convex 
polyhedral cone and let $\square\in \cC(\sigma^\vee)$.
Assume that there exists a log discrepancy function 
$\psi\colon \sigma \to \R$ such that $\square$ is
asymptotically $\psi$-saturated.

Then for every $e_1\in \sigma\setminus 0$, 
there exist $m\in M_\Q\cap \square$ 
and a rational convex polyhedral cone 
$\sigma_1\subset \sigma$, with the following 
properties:
\begin{itemize}
\item[(i)] $e_1\in \relint(\sigma_1)$.
\item[(ii)] $h_\square(e)=
\langle m,e\rangle$ for $e\in \sigma_1$.
\end{itemize}
\end{thm}

\begin{proof} Choose norms $\Vert\cdot \Vert$ on $N_\R$ 
and $M_\R$, defined as the maximum of the absolute values 
of the components with respect to some basis of $N$ and 
its dual basis in $M$, respectively. Let 
$$
S(\sigma)=\{e\in \sigma;\Vert e\Vert=1\}.
$$
Define the positive real number $\epsilon(\psi)$ by the 
formula
$$
-\epsilon(\psi)^{-1}=\min\{\frac{\langle m,e\rangle}{\psi(e)};
\Vert m\Vert =1, e\in S(\sigma) \}.
$$
The restriction of $\psi$ to $S(\sigma)$ is a positive, 
continuous function, hence $\epsilon(\psi)$ is a well defined. 
In particular, 
$$
\langle m,e\rangle+\psi(e)>0 \mbox{ for } 
0\ne e\in \sigma , \Vert m\Vert<\epsilon(\psi).
$$	
Denote by $h\in \cS(\sigma)$ the support function of $\square$.
(1) There exists $m\in M_\Q\cap \square_h$ such that 
$\langle m,e_1\rangle=h(e_1)$.
\newline

Indeed, let $\tau$ be the unique face of $\sigma$ 
which contains $e_1$ in its relative interior. We may 
find orthogonal decompositions 
$$
N=N'\oplus N'', M=M'\oplus M'',
$$ 
where $N'=N\cap(\tau-\tau), M'=M\cap \tau^\perp$
and $M',N'$ and $M'',N''$ are dual lattices, respectively.
If $N''\ne 0$, let $\sigma''$ be the image of 
$\sigma$ under the projection map $N_\R\to N''_\R$. 
Since $\tau\supseteq \sigma\cap(-\sigma)$, we infer that
$\sigma''$ is a strongly rational convex polyhedral cone. 

(1a) Since $h$ is the support function of the 
non-empty convex set $\square_h$, there exists a 
sequence of points $m_k\in \square_h$ such that
$$
\lim_{k\to \infty}\langle m_k,e_1\rangle=h(e_1).
$$
If we decompose $m_k=m'_k+m''_k$, we claim that 
$m'_k$ belongs to a bounded set of $M'_\R$. 

Indeed, assume by contradiction that 
$\lim_{k\to \infty}\Vert m'_k\Vert=+\infty$. 
By the usual compactness argument, we may assume that 
there exists $m'\in M'_\R$ such that
$$
\lim_{k\to \infty} \frac{m'_k}{\Vert m'_k\Vert}=m'.
$$
For every $e\in \tau$, we have
$$
\langle \frac{m'_k}{\Vert m'_k\Vert},e\rangle=
\langle \frac{m_k}{\Vert m'_k\Vert},e\rangle\ge
\frac{h(e)}{\Vert m'_k\Vert}.
$$
Letting $k$ converge to infinity, we obtain
$\langle m',e\rangle\ge 0$. Furthermore,
$$
\lim_{k\to \infty}\langle m'_k,e_1\rangle=
\lim_{k\to \infty}\langle m_k,e_1\rangle=h(e_1),
$$
so a similar argument gives $\langle m',e_1\rangle=0$.
Therefore $0\ne m'\in \tau^\vee\cap e_1^\perp$.
Since $e_1$ belongs to the relative interior of $\tau$,
we infer that $m'\in \tau^\perp$. This implies
$m'=0$, a contradiction. Therefore the claim holds.

(1b) By (1a), we may replace $(m_k)_k$ by a subsequence
so that there exists $m'\in M'_\R$ with
$$
\lim_{k\to \infty} m'_k=m'\mbox{ and }
\langle m',e_1\rangle=h(e_1).
$$
By Theorem~\ref{appusual}, there exists a positive 
multiple $j$ of $I$ such the following system has a 
solution:
\[
\left\{ \begin{array}{l}
 m'_j\in M' \\
 \langle jm',e_1\rangle-\psi(e_1)<\langle m'_j,e_1\rangle 
 \le \langle jm',e_1\rangle \\
 \Vert m'_j-j m'\Vert<\frac{1}{2}\epsilon(\psi)
\end{array}
\right. \]
We now choose $k$ large enough so that
$$
j\Vert m'-m'_k\Vert<\frac{1}{2}\epsilon(\psi).
$$
Since $\sigma''$ is a strongly rational convex
polyhedral cone, the following system has a solution
\[
\left\{ \begin{array}{l}
 m''_j\in M'' \\
 m''_j\in jm''_k+{\sigma''}^\vee
\end{array}
\right. \]
Set $m_j=m'_j+m''_j\in M$. The following holds for $e\in S(\sigma)$:
\begin{align*}
\langle m_j,e\rangle-j h(e)+\psi(e)  & =
\langle m_j-j m_k,e\rangle+\psi(e)+j(\langle m_k,e\rangle-h(e))\\
   & \ge \langle m_j-j m_k,e\rangle+\psi(e)\\
   & \ge \langle m'_j-j m'_k,e\rangle+\psi(e)\\
   & > 0,
\end{align*}
where the latter inequality follows from
$$
\Vert m'_j-j m'_k\Vert\le \Vert m'_j-j m'\Vert+ 
\Vert jm'-j m'_k\Vert<\epsilon(\psi).
$$
Since $jh$ is $\psi$-saturated, we infer that
$m_j\in \square_{jh}$. In particular,
$$
\langle m_j,e_1\rangle\ge j h(e_1).
$$
The opposite inequality holds from construction, hence 
$
\langle m_j,e_1\rangle=j h(e_1).
$ 
Therefore we obtain
$$
m:=\frac{m_j}{j}\in M_\Q\cap \square_h,
\langle m,e_1\rangle=h(e_1).
$$

(2) Since $m$ is rational, we may replace $\square$ by
$\square-m$, or equivalently, we replace $h$ by $h-m$.
Thus we may assume that $0\in \square_h$ and
$e_1\in \sigma_0$, where
$$
\sigma_0=\{e\in \sigma; h(e)=0\}.
$$  
There exists a decomposition $N=N'\oplus N''$, with
dual decomposition $M=M'\oplus M''$, such that
$e_1=e'_1+e''_1$, $e'_1\in N'_\Q, e''_1\in N''_\R$ and 
$$
\{m''\in M''; \langle m'',e''_1\rangle\in \Q\}=\{0\}.
$$
If $e''_1=0$, then $e_1\in N_\Q$ and the theorem holds for
$\sigma_1=\R_{\ge 0}\cdot e_1$ and $m=0$. 

(2a) Assume that $e''_1\ne 0$. We claim that the 
following equality holds
$$
\sigma_0^\vee\cap e_1^\perp=\sigma_0^\vee\cap 
M'_\R\cap {e'_1}^\perp.
$$
We only have to prove the direct inclusion. Fix
$m\in \sigma_0^\vee\cap e_1^\perp$. 
We have to show that $m''=0$, where $m=m'+m''$ is the
decomposition in $M'_\R\oplus M''_\R$. Assume by
contradiction that $m''\ne 0$. Since
$\langle m,e_1\rangle\in \Q$, we infer by 
Lemma~\ref{app2} and Theorem~\ref{app1} that there 
exist a positive integer $k$ and $m_1\in M$ such that 
$-\psi(e_1)<\langle m_1-km,e_1\rangle<0$ and
$\Vert m_1-km\Vert<\epsilon(\psi)$. Since
$\langle m,e_1\rangle=0$, we obtain
$$
-\psi(e_1)<\langle m_1,e_1\rangle<0,
\Vert m_1-km\Vert<\epsilon(\psi).
$$
We consider the following set 
$$
C=\{e\in S(\sigma); \langle m_1,e\rangle+\psi(e)\le 0\}.
$$
Since $0\in \square_h$, we have $h\le 0$. 
If the set $C$ is empty, then 
$$
\langle m_1,e\rangle-jh(e)+\psi(e)\ge 
\langle m_1,e\rangle+\psi(e)>0, \forall e\in S(\sigma).
$$
Therefore $m_1\in M\cap\stackrel{\circ}{\square}_{jh-\psi}$, 
hence $m_1\in \square_{jh}$ by saturation. In particular, 
$\langle m_1,e_1\rangle\ge 0$, which contradicts the choice
of $m_1$. 

Therefore the set $C$ is non-empty. Since $\psi$
is continuous, $C$ is also compact. If 
$C\cap \sigma_0=\emptyset$, then there exists a positive
integer $j$ with 
$$
j>\sup_{e\in C}\frac{\langle m_1,e\rangle+\psi(e)}{h(e)}.
$$
Then $m_1\in M\cap \stackrel{\circ}{\square}_{jh-\psi}$,
and saturation implies that $m_1\in \square_{jh}$,
hence $\langle m_1,e_1\rangle\ge 0$, which contradicts the 
choice of $m_1$.

Therefore there exists $e\in C\cap \sigma_0$. In particular,
$$
\langle m,e\rangle < \frac{\langle m_1,e\rangle+\psi(e)}{k}
\le 0.
$$
Therefore $\langle m,e\rangle<0$, contradicting the assumption 
$e\in \sigma_0, m\in \sigma_0^\vee$.

(2b) The function $h$ is continuous, being upper convex
(\cite{Rock70}, Theorem 10.1).
Therefore $\sigma_0$ is a closed convex cone in $N_\R$. 
By duality (cf.~\cite{Oda88}, Theorem A.1), (2a) is equivalent 
to
$$
\sigma_0+\R\cdot e_1'+N''_\R=\sigma_0+\R\cdot e_1.
$$
In particular, there exists an open neighborhood $U''$ of 
$e_1''$ in $N''_\R$ such that
$$
e_1'+U''\subset \sigma_0.
$$
Since $\dim(U'')=\dim(N'')$, there exist 
$\bar{e}_1,\ldots,\bar{e}_{n+1}\in U''\cap N''_\Q$, 
where $n=\dim(N'')$, and there exists $\lambda_i\in (0,1)$ 
such that $\sum_{i=1}^{n+1}\lambda_i=1$ and 
$$
e''_1=\sum_{i=1}^{n+1}\lambda_i \bar{e}_i.
$$
Let $\sigma_1$ be the rational polyhedral cone spanned
by $e'_1+\bar{e}_1,\ldots,e'_1+\bar{e}_{n+1}$. It is clear that
$\sigma_1\subset \sigma$, 
$e_1\in \relint(\sigma_1)$ and $h\vert_{\sigma_1}=0$.
\end{proof}

\begin{thm}\label{mainfga} 
Let $\sigma \subseteq N_\R$ be a rational convex 
polyhedral cone and let $\square\in \cC(\sigma^\vee)$, with
support function $h\in \cS(\sigma)$.
Assume that there exists a log discrepancy function
$\psi\colon \sigma \to \R$ such that 
\begin{itemize}
\item[(i)] $\square$ is asymptotically $\psi$-saturated;
\item[(ii)] $h-\psi$ is upper convex.
\end{itemize}
Then $\square$ is a rational convex polyhedral set.
\end{thm}

\begin{proof} 
We prove the result by induction on $\dim(N)$. 
If $\dim(N)=1$, then $\square$ is either a point
or an interval of the form $[a,b]$ or $[a,+\infty)$. 
Its endpoints are rational by Theorem~\ref{1rational}, 
hence $\square_h$ is a rational convex polyhedral set.

Assume now that $\dim(N)>1$ and the theorem holds for
smaller dimensional lattices $N$. 
We prove the theorem in three steps.

\begin{itemize}
\item[(1)] Assume $0\in \square$ and 
$\sigma_1\subset \sigma$ is a rational
convex polyhedral cone such that $h\vert_{\sigma_1}=0$.
Then there exists a rational polyhedral cone 
$\sigma_2\subset \sigma$ such that
$\relint(\sigma_1)\subset \relint(\sigma_2)$, and one of 
the following two properties holds:
\begin{itemize}
 
\item[(a)] $\dim(\sigma_2)=\dim(\sigma_1)+1$, and
$h\vert_{\sigma_2}=0$.

\item[(b)] $\dim(\sigma_2)=\dim(\sigma)$ and there 
exist finitely many rational points
$m_1,\ldots, m_n\in M_\Q\cap \square$ such that
for every $e\in \sigma_2$ there exists some $i$
with $\langle m_i,e\rangle=h(e)$.
\end{itemize}

\begin{proof} Let $N'=N/(N\cap(\sigma_1-\sigma_1))$, with
dual lattice $M'=M\cap \sigma_1^\perp$. 
If $\dim(\sigma_1)=\dim(N)$, we are in case (1b). Assume 
now $\dim(\sigma_1)<\dim(N)$, so that $0<\dim(N')<\dim(N)$.
With the notations of Lemma~\ref{restriction}, we 
have a projection homomorphism
$$
\pi_\Z\colon N\to N', \sigma'=\pi(\sigma),
$$
the support function $h'\colon \sigma'\to \R$ of 
$\square\cap M'_\R$ and the log discrepancy functions 
$\psi'_k \colon \sigma'\to \R $, for $k\ge 1$. 
By Lemma~\ref{restriction}, $\square_{h'}$ is asymptotically 
$\psi'_k$-saturated and $kh'-\psi'_k$ is upper convex.
The inductive assumption implies that there exists a finite 
set $\{m'_i\}_{i\in I}\subset M'_\Q\cap \square$ such that 
for every $e'\in \sigma'$, 
$h'(e')=\langle m'_i,e'\rangle$ for some $i\in I$.
We distinguish two cases, depending on whether the convex
set $\square_{h'}\subseteq M'_\R$ is maximal dimensional 
or not.

(a) Assume $\dim(\square_{h'})<\dim(M')$. Equivalently, the
lattice $N''=N'\cap {\square_{h'}}^\perp$ is non-zero. 
Let $\psi''_k=\psi'_k\vert_{N''}$ and let $M''$ be the dual 
lattice of $N''$. By Theorem~\ref{chr}, 
$M''\cap \stackrel{\circ}{\square}_{-\psi''_k}=\{0\}$.
Furthermore, $kh-\psi_k$ is upper convex and $h\vert_{N''_\R}$ 
is linear, hence $-\psi''_k$ is upper convex. Therefore
Theorem~\ref{klbound} applies, hence there exists 
$0\ne e'_k\in N''$ such that
$$
\psi''_k(e'_k)+\psi''_k(-e'_k)\le C,
$$ 
where $C$ is a positive constant depending only on 
$\dim(N'')$. By Lemma~\ref{4th}, the $e'_k$'s belong to a 
compact set, hence we may assume that $e'_k=e'$ for infinitely
many $k$'s. Then there exist $e^+_k,e^-_k\in \sigma$ such that
$\pi(e^+_k)=e'$, $\pi(e^-_k)=-e'$ and
$$
k[h'(e')-h(e^+_k)+h'(-e')-h(e^-_k)]+
\psi(e^+_k)+\psi(e^-_k)\le C+1.
$$
In particular,
$$
\psi(e^+_k)+\psi(e^-_k)\le C+1.
$$
By Lemma~\ref{4th}, the sequences $(e_k^+)_k, (e_k^-)_k$
belong to a compact set, so we may assume that the limits
$e^-=\lim_{k\to \infty} e_k^-, e^+=\lim_{k\to \infty} e^+_k$
exist. It is clear that $e^+,e^-\in \sigma$
and $\pi(e^+)=e',\pi(e^-)=-e'$.
The above inequality and the positivity of $\psi$
implies
$$
h'(e')-h(e^+_k)+h'(-e')-h(e^-_k)\le \frac{C+1}{k}.
$$
Letting $k$ converge to infinity, we obtain 
$h'(e')=h(e^+)$, $h'(-e')=h(e^-)$. Since $e'\in N''$, we 
have $h'(e')=h'(-e')=0$. 

We claim that we may assume that $e^+,e^-\in N_\Q$. Indeed,
since $\pi(e^+)\in N'_\Q$ and $\sigma_1$ is rational, there
exists $f\in \sigma_1$ such that $e^++f\in \sigma_1\cap N$.
Then $h(e^++f)\ge h(e^+)+h(f)=0$, hence $h(e^++f)=0$.
Also, $\pi(e^++f)=e'$, hence we may replace $e^+$ by
$e^++f$. A similar argument applies to $e^-$.

It is easy to verify that the rational convex polyhedral 
cone $$\sigma_2=\sigma_1+\R_{\ge 0}e^++\R_{\ge 0}e^-$$
satisfies (1a).

(b) Assume $\dim(\square_{h'})=\dim(M')$. In this case, 
the ample fan $\Delta_{h'}$ of $h'$ is a fan in $N'$ with 
$\vert \Delta_{h'}\vert=\sigma'$.

(b1) For every $e'\in \Delta_{h'}(1)$, there exists 
$e\in \sigma\cap N$ such that $\pi(e)=e'$ and $h(e)=h'(e')$.

Indeed, since $h'$ is rational piecewise linear and
asymptotically $\psi'_k$-saturated, we obtain by 
Theorem~\ref{chr} that 
$$
\psi'_k(e')\le 1 , \forall k\ge 1.
$$ 
Therefore there exists $e_k\in \sigma$ such that
$\pi(e_k)=e'$ and 
$$
kh'(e')-(kh-\psi)(e_k)\le 2.
$$
In particular, we obtain
$
\psi(e_k)\le 2.
$
By Lemma~\ref{4th}, the sequence $(e_k)_k$ belongs to a 
bounded set, so we may assume that the limit 
$e=\lim_{k\to \infty} e_k$ exists. We clearly have 
$\pi(e)=e'$. The positivity of $\psi$ implies
$$
h'(e')-h(e_k)\le \frac{2}{k}.
$$
Letting $k$ converge to infinity, we obtain
$
h'(e')-h(e)\le 0,
$
hence $h'(e')-h(e)=0$. The rationality of $e$ is obtained
the same way as in the proof of (a) above.

(b2) Let $\tau$ be a maximal dimensional cone of $\sigma'$,
spanned by $e'_1,\ldots,e'_r\in \Delta_{h'}(1)$. There
exists $i\in I$ such that $h'(e')=\langle m_i,e'\rangle$
for every $e'\in \sigma'$. By (b1), there exist 
$e_j\in \sigma\cap N_\Q$ such that
$\pi(e_j)=e'_j$ and $h(e_j)=h'(e'_j)$, for $1\le j\le r$. 
Therefore 
$
h(e)=\langle  m_i,e'\rangle
$ 
for every 
$
e\in \sigma_1+\sum_{j=1}^p \R_{\ge 0}e_i.
$
The cone $\sigma_1+\sum_{j=1}^p \R_{\ge 0}e_i\subset \sigma$
has the same dimension as $\sigma$. The union of all these 
cones, taken after all maximal cones $\tau$ in $\Delta$,
contains a cone $\sigma_2$ satisfying (1b) with respect to
$\{m'_i\}_{i\in I}$.
\end{proof}

\item[(2)] Every non-zero point $e\in \sigma$ has an open 
polyhedral neighborhood on which $h$ is rational,
piecewise linear.

Indeed, fix $e$ as above. By Theorem~\ref{1rational},
there exists $m_0\in M_\Q\cap \square_h$ and
there exists a rational convex polyhedral cone
$\sigma_0\subset \sigma$ such
that $e\in \relint(\sigma_0)$ and $h(e)=\langle 
m_0,e\rangle$ for every $e\in \sigma_0$. 

We may replace $\square_h$ by its rational translate
$
\square_h-m_0,
$
so that we may assume that $m_0=0$. In particular,
$0\in \square_h$ and $h\vert_{\sigma_0}=0$. By (1), 
either the claim holds, or there exists a 
$(\dim(\sigma_0)+1)$-dimensional rational polyhedral 
cone $\sigma_1\subseteq\sigma$ such that 
$\relint(\sigma_0)\subset \relint(\sigma_1)$ and
$h\vert_{\sigma_1}=0$.
By (1) again, either the claim holds, or there
exists a $(\dim(\sigma_1)+2)$-dimensional cone 
$\sigma_2\subseteq \sigma$ such that 
$\relint(\sigma_1)\subset \relint(\sigma_2)$
and $h\vert{\sigma_2}=0$.
We repeat this procedure for $\sigma_2$ and so on. This 
procedure clearly stops in a finite number of steps, 
hence the claim holds.

\item[(3)] Fix a norm $\Vert\cdot \Vert$ on $N_\R$ and
set
$
S(\sigma)=\{e\in\sigma;\Vert e\Vert=1\}.
$ 
For each point $e\in S(\sigma)$, we consider the pair 
$(\sigma(e); \{m_i(e)\}_{i\in I(e)})$ 
constructed in (2). We obtain an open covering
$$
S(\sigma)=\bigcup_{e\in S(\sigma)}
S(\sigma)\cap\relint(\sigma(e)).
$$
Since $S(\sigma)$ is compact,
it may be covered by the relative interiors
of the cones corresponding to finitely many points
$e_1,\ldots, e_k$. Let $K$ be the convex hull of the 
finitely many rational points
$$
\{m_i(e_1)\}_{i\in I(e_1)}\cup \ldots 
\cup\{m_i(e_k)\}_{i\in I(e_k)}.
$$
Then 
$
\square=K+\sigma^\vee,
$ 
i.e. $\square$ is a rational convex polyhedral set.
\end{itemize}
\end{proof}


\section{Toric FGA algebras}


\begin{thm}\label{maini}
Let $\pi\colon X\to S$ be a proper surjective toric 
morphism with connected fibers, and let $B$ be an 
invariant $\Q$-divisor on $X$ such that $(X,B)$ is a 
log pair with Kawamata log terminal singularities 
and $-(K+B)$ is $\pi$-nef. Let 
$$
\cL\subseteq \bigoplus_{i=0}^\infty\pi_*\cO_X(iD)
$$ 
be an invariant graded $\cO_S$-subalgebra, where $D$ 
is an invariant $\R$-divisor on $X$, such that $\cL$ 
is asymptotically saturated with respect to $(X/S,B)$.

Then $\cL$ is finitely generated.
\end{thm}

\begin{proof} We may assume that $S$ is affine.
Then $X=T_N\emb(\Delta),S=T_{\bar{N}}(\bar{\sigma})$,
and $\pi$ corresponds to a map of fans 
$\varphi_\Z\colon (N,\Delta)\to (\bar{N},\bar{\sigma})$
such that $\vert\Delta\vert=\varphi^{-1}(\bar{\sigma})$
is a rational convex set, denoted by $\sigma$.

We can write $B=\sum_{e\in \Delta(1)}b_e V(e)$, where
$\Delta(1)$ is the set of primitve vectors on the one 
dimensional cones of $\Delta$. The log canonical divisor
$K+B$ is represented by a function $\psi\colon \sigma\to 
\R$ such that $\psi$ is $\Delta$-linear and $\psi(e)=1-b_e$
for every $e\in \Delta(1)$. Since $(X,B)$ has Kawamata log
terminal singularities, $\psi$ is a log discrepancy function.

If $\cL=\cL_0$, then $\cL$ is finitely generated. Otherwise,
we may replace $I$ by a multiple so that $\cL_i\ne 0$ for 
every $I\vert i$.
Let $i$ be a positive multiple of $I$. Since $\cL$ is 
torus invariant, there exist finitely many lattice points
$m_{i,1},\ldots,m_{i,n_i}\in M$ such that $\chi^{m_{i,1}},
\ldots, \chi^{m_{i,n_i}}$ generate the $\cO_S$-module $\cL_i$.
Define $h_i\colon \sigma\to \R$ by 
$$
h_i(e)=\min_{j=1}^n \langle m_{i,j},e\rangle.
$$ 
The support function $h_i$ is independent of the choice 
of generators, and the torus invariant $\cO_S$-algebra
$$
\bar{\cL}=\bigoplus_{i=0}^\infty 
(\bigoplus_{m\in M\cap \square_{h_i}}{\mathbb C} \cdot \chi^m)
$$
is the integral closure of $\cL$ in its field of fractions 
(\cite{PLflips}, Proposition 4.15).

Choose a refinement $\Delta_i$ of the 
fan $\Delta$ so that $\Delta_i$ is a simple fan and $h_i$ 
is $\Delta_i$-linear. This corresponds to a toric resolution 
of singularities $\mu_i\colon X_i=T_N\emb(\Delta_i)\to X$ such that 
$M_i=\sum_{e\in \Delta_i(1)}-h_i(e)V(e)$ is a $\pi\circ \mu_i$-free 
divisor. Since $X_i$ is nonsingular, the union of its invariant
prime divisors $\sum_{e\in \Delta_i(1)}V(e)$ has simple normal 
crosssings.

In the above set-up, $\cL$ is asymptotically saturated with 
respect to $(X/S,B)$ if and only if 
$$
H^0(X_i,\lceil K_{X_i}-\mu_i^*(K+B)+\frac{j}{i}M_i\rceil)
\subseteq H^0(X_j,M_j), \forall I\vert i,j.
$$

{\em Step 1}: Asymptotic saturation is equivalent to the 
following property:
$$
M\cap \stackrel{\circ}{\square}_{\frac{j}{i}h_i-\psi}
\subset \square_{h_j}, \forall I\vert i,j.
$$
Indeed, let $m\in M$. Then
$\chi^m\in H^0(X_i,\lceil K_{X_i}-\mu_i^*(K+B)+
\frac{j}{i}M_i\rceil)$ if and only if 
$\langle m,e\rangle+\lceil -1+\psi(e)-\frac{j}{i}h_i(e)
\rceil\ge 0$ for every $e\in \Delta_i(1)$. Since 
$\langle m,e\rangle\in \Z$, this is equivalent to
$
\langle m,e\rangle>\frac{j}{i}h_i(e)-\psi(e)
$
for every $e\in \Delta_i(1)$. Since $\psi$ and $h_i$ are
$\Delta_i$-linear, the latter is equivalent to 
$
\langle m,e\rangle>\frac{j}{i}h_i(e)-\psi(e)
$
for every $e\in \sigma\setminus 0$. 
On the other hand, $\chi^m\in H^0(X_j,M_j)$ if and only 
if $m\in \square_{h_j}$. 

{\em Step 2}: The function $h=\lim_{i\to \infty}
\frac{1}{i}h_i\colon \sigma\to \R$ is a well defined
positively homogeneous, upper convex function.

Indeed, we can write $D=\sum_{e\in \Delta(1)}d_eV(e)$.
Let $\tilde{h}\colon \sigma\to \R$ be the support function
of the convex set 
$$
\{m\in M_\R; \langle m,e\rangle\ge -d_e,\forall
e\in\Delta(1)\}.
$$
Since $\cL_i\subseteq H^0(X,iD)$, we obtain
$h_i\ge i\tilde{h}$. On the other hand, the property
$\cL_i\cdot \cL_j\subseteq \cL_{i+j}$ implies 
$h_i+h_j\ge h_{i+j}$. Then it is easy to see that for
every $e\in \sigma$, the sequence $\frac{1}{i}h_{i}(e)$ 
is bounded from below and converges to its infimum. 
Being a limit of positively homogeneous
upper convex functions, $h$ satisfies these two properties
too. Note that $h_i\ge ih$ for every $i$.

{\em Step 3}:  Asymptotic saturation is equivalent to 
the following property:
$$
M\cap \stackrel{\circ}{\square}_{j h-\psi}
\subset \square_{h_j}, \forall I\vert j.
$$
Indeed, fix $I\vert j$, choose a norm $\Vert \cdot \Vert$ 
on $N_\R$ and set $S(\sigma)=\{e\in\sigma; \Vert e\Vert=1\}$.
Let $m\in M \cap \stackrel{\circ}{\square}_{j h-\psi}$.
This means that the function 
$$
f\colon S(\sigma)\to \R, e\mapsto 
\langle m,e\rangle-j h(e)+\psi(e)
$$
is positive. The functions $\frac{1}{i}h_i$ are upper
convex, hence they converge uniformly to $h$ on the 
compact set $S(\sigma)$, by Theorem~\ref{Rocka}.
Therefore there exists some $i$ such that
$
f-j(\frac{1}{i}h_i-h)\vert_{S(\sigma)}
$
is a positive function. This means that 
$m\in M\cap \stackrel{\circ}{\square}_{\frac{j}{i}h-\psi}$. 
By Step 1, we obtain $m\in \square_{h_j}$.

The converse is clear by Step 1, since $h_i\ge ih$.

{\em Step 4}: $nh=h_n$ for some $I\vert n$.

Indeed, $-\psi$ is upper convex since $-(K+B)$ is nef. 
Therefore $h-\psi$ is upper convex. By Step 3, the
hypothesis of Theorem~\ref{mainfga} are satisfied,
hence $\square_h$ is a rational polyhedral set. In 
particular, there exists some $I\vert n$ such that 
$\square_{nh}$ is the convex hull of its lattice points.
We have
$$
M\cap \square_{nh}\subseteq 
M\cap \stackrel{\circ}{\square}_{nh-\psi}
\subset \square_{h_n},
$$
hence $\square_{nh}$, the convex hull of $M\cap \square_{nh}$, 
is included in $\square_{h_n}$. Therefore $h_n\ge nh$. The
opposite inclusion holds by construction, hence $nh=h_n$.

{\em Step 5}:
We have $kh_n\ge h_{kn}\ge knh$ for $k\ge 1$. By Step 4,
we obtain $h_{kn}=kh_n$ for every $k\ge 1$. This means that
$$
\bigoplus_{k=0}^\infty\bar{\cL}_{kn}=\bigoplus_{k=0}^\infty
(\pi\circ \mu_n)_*\cO_{X_n}(kM_n).
$$
The right hand side is finitely generated since $M_n$ is a
$\pi\circ \mu_n$-free divisor, hence 
$\bigoplus_{k=0}^\infty\bar{\cL}_{kn}$ is also finitely generated.
Therefore $\bar{\cL}$ is finitely generated (cf.~\cite{PLflips},
Theorem 4.6). Since $\bar{\cL}$ is the integral 
closure of $\cL$ in its field of fractions, we conclude that
$\cL$ is finitely generated.
\end{proof}

\begin{proof}(of Theorem~\ref{main}) The statement is local
over $S$, hence we may assume that $S$ is affine. 
Thus $X=T_N\emb(\Delta),S=T_{\bar{N}}(\bar{\sigma})$,
and $\pi$ corresponds to a map of fans 
$\varphi_\Z\colon (N,\Delta)\to (\bar{N},\bar{\sigma})$
such that $\vert\Delta\vert=\varphi^{-1}(\bar{\sigma})$
is a rational convex set.

Let $\cL\subseteq \bigoplus_{i=0}^\infty\pi_*\cO_X(iD)$ 
be an invariant graded $\cO_S$-subalgebra, where $D$ 
is an invariant $\R$-divisor on $X$, such that $\cL$ 
is asymptotically saturated with respect to $(X/S,B)$.
Then $\cL$ is finitely generated by Theorem~\ref{maini},
which proves (1).
In particular, there exists a rational convex polyhedral set
$\square\in \cC(\vert\Delta\vert^\vee)$ (corresponding
to the limit support function $h\in \cS(\vert\Delta\vert)$ 
in the Step 2 of the proof of Theorem~\ref{maini})
and an $S$-isomorphism
$$
\Proj(\bar{\cL})\simeq \Proj(\bigoplus_{i=0}^\infty
\bigoplus_{m\in M\cap i\square}{\mathbb C}\cdot \chi^m).
$$
The right hand side is the torus embedding of the ample
fan $\Delta_\square$. Let $\psi\colon \vert\Delta\vert
\to \R$ be the log discrepancy function associated to
$(X,B)$ (Example~\ref{gld}). Then $\square$ is 
asymptotically $\psi$-saturated. Since $-(K+B)$ is 
$\pi$-nef, $-\psi$ is upper convex. Therefore 
$h_\square-\psi$ is upper convex. Theorem~\ref{asyccs}
applies, hence $\Delta_\square$ belongs to a finite
set of fans associated to $(X/S,B)$. This proves (2).

\end{proof}


\end{document}